\documentclass[reqno,11pt,a4paper]{amsart}
\usepackage{graphicx,color,amssymb,latexsym,amsfonts,amsmath}
\usepackage{mathrsfs}
\usepackage{graphics}
\usepackage{lscape}
\usepackage{dsfont}
\usepackage[normalem]{ulem} 
\usepackage{verbatim}

\renewcommand{\b}[1]{\mbox{\boldmath $#1$}}

\newtheorem{proposition}{\bf Proposition }[section]%

\newtheorem{theorem}{\bf Theorem}[section]

\numberwithin{equation}{section} \theoremstyle{plain}
\theoremstyle{definition}

\newtheorem{remark}{Remark}

\newcommand{\E}{\mathsf{E}}
\renewcommand{\P}{\mathsf{P}}

\begin{document}

\title[] { Asymptotic analysis of ruin in CEV model}
\author{F. Klebaner}
\address{School of Mathematical Sciences,\\ Building 28M, Monash
University, \\ Clayton Campus, Victoria 3800,\\ Australia}
\email{fima.klebaner@sci.monash.edu.au}
\thanks{}

\author{R. Liptser}
\address{Department of Electrical Engineering Systems,
Tel Aviv University, 69978 Tel Aviv, Israel}
\email{liptser@eng.tau.ac.il; rliptser@gmail.com}


\maketitle

\begin{abstract}
We give asymptotic analysis for probability of absorbtion
$\mathsf{P}(\tau_0\le T)$ on the interval $[0,T]$, where $
\tau_0=\inf\{t:X_t=0\}$ and $X_t$ is a nonnegative diffusion
process relative to Brownian motion $B_t$,
\begin{align*}
dX_t&=\mu X_tdt+\sigma X^\gamma_tdB_t.
\\
X_0&=K>0
\end{align*}
Diffusion parameter $\sigma x^\gamma$, $\gamma\in [\frac{1}{2},1)$
is not Lipschitz continuous and assures $\mathsf{P}(\tau_0>T)>0$.
Our main result:
$$
\lim\limits_{K\to\infty}
\frac{1}{K^{2(1-\gamma)}}\log\mathsf{P}(\tau_{0 }\le T)
=-\frac{1}{2\E M^2_T},
$$
where $ M_T=\int_0^T\sigma(1-\gamma)e^{-(1-\gamma)\mu s}dB_s $.
Moreover we describe the most likely path to absorbtion of the
normed process $\frac{X_t}{K}$ for $K\to\infty$.

\end{abstract}

\maketitle
\section{\bf Introduction}
\label{sec-1}
 In this paper, we analyze the
\textit{Constant Elasticity of Variance Model} (CEV), introduced
by Cox 1996, \cite{Cox} and applied to \textit{Option Pricing}
(see e.g. Delbaen, Shirakawa \cite{DelShir} and Lu, Hsu
\cite{LuChia}). This model is given by the It\^o equation with
respect  to a standard Brownian motion $B_t$ and a positive
initial condition $X_0=K>0$,
\begin{equation}\label{eq:1x}
dX_t=\mu X_tdt+\sigma X^\gamma_tdB_t,
\end{equation}
where $\mu,\sigma\ne 0$ are arbitrary constants and $\gamma \in
[\frac{1}{2},\ 1)$. For $\gamma=\frac{1}{2}$, this model is known
as CIR (Cox, Ingersol and Ross) model.  The diffusion coefficient
$\sigma x^\gamma$ is only H\"older continuous, yet the It\^o
equation \eqref{eq:1x} has   unique strong
solution\footnote{Delbaen and Shirakawa, \cite{DelShir} -
existence; Yamada-Watanabe - uniqueness  (see e.g., Rogers and
Williams, p. 265 \cite{RW} or \cite{Kleb} p.17 and Theorem 13.1)}.
In contrast to Black-Scholes model ($\gamma=1$) with $X_t>0$ for
any $t>0$, for CEV model the process $X_t$ is absorbed in zero at
the time $ \tau_0=\inf\{t:X_t=0\} $ with  $\P(\tau_0<\infty)>0$
which can be interpreted as time of ruin.

In a proposed asymptotic analysis, as $K\to\infty$, a crucial role
plays the normed process $ x^K_t=\frac{X_t}{K}, $ being the unique
solution of the following It\^o equation
\begin{equation}\label{eq:xKt}
dx^K_t=\mu x^K_tdt+\frac{\sigma}{K^{(1-\gamma)}}(x^K_t)^\gamma
dB_t,
\end{equation}
subject to the initial condition $x^K_0=1$, and a small diffusion
parameter $\frac{\sigma x^\gamma}{K^{1-\gamma}}$. We emphasize
that the process $x^K_t$ inherits the ruin time $\tau_0$.

The assumption  $\gamma<1$ implies that the diffusion in
\eqref{eq:xKt} has a small diffusion coefficient. This enables us
to find a rough lower bound of $\P(\tau_0\le T)$ for any $K>0$
(see Remark \ref{rem-1}). With  $K\to\infty$, this lower bound
 is best possible on  logarithmic scale.
To this end we apply the Large Deviation Theory for asymptotic
analysis of two families:
$$
\Big\{\Big(x^K_t\Big)_{t\in[0,T]}\Big\}_{K\to\infty}
\quad\text{and}\quad
\Big\{\frac{1}{K^{1-\gamma}}M_T\Big\}_{K\to\infty},
$$
where
\begin{equation}\label{Meq}
M_t=\int_0^t\sigma(1-\gamma)e^{-(1-\gamma)\mu s}dB_s.
\end{equation}
For the second family, Large Deviation Principle (LDP) is well
known. For the fist family, Freidlin-Wentzell's LDP, \cite{FW}, is
anticipated even though the diffusion parameter is only H\"older
continuous and  singular at zero. For $\gamma=\frac{1}{2}$, LDP is
known from Donati-Martin et al.,\cite{Yor}, with the speed rate
$\frac{1}{K}$ and the rate function of Freidlin-Wentzell's type
with a corresponding modification:
$
J_T(u)=\frac{1}{2}\int_0^T\big(\frac{\dot{u}_t-\mu
u_t}{\sqrt{u_t}}\big)^2I_{\{u_t>0\}}dt. $ \footnote{For
$\gamma=\frac{1}{2}$ see also  \cite{KLIP} and Rouault \cite{Rou}.
} We show that for $\gamma\in \big(\frac{1}{2},1\big)$ LDP is also
valid with the speed rate and the rate function depending on
$\gamma$. Combining both LDP's  we obtain the following asymptotic
result: there is a smooth nonnegative function $u^*_t$, with
$u^*_0=1$ and absorbed at the time   $T$, $u^*_T=0$, such that for
any smooth nonnegative function $u_t$, with $u_0=1$ and absorbed
on the interval $[0,T]$, $u_T=0$,
\begin{gather}
\lim_{K\to\infty}\frac{1}{K^{2(1-\gamma)}}\log\P\big(\tau_0\le T\big)
\overrightarrow{}\nonumber\\
=\lim_{\delta\to 0}\lim_{K\to\infty}\frac{1}{K^{2(1-\gamma)}}\log\P\bigg(\sup_{t\in[0,T]}
\big|x^K_t-u^*_t\big|\le \delta
\bigg)
\nonumber\\
\ge \lim_{\delta\to 0}\lim_{K\to\infty}\frac{1}{K^{2(1-\gamma)}}\log\P\bigg(\sup_{t\in[0,T]}
\big|x^K_t-u_t\big|\le \delta
\bigg).
\nonumber
\end{gather}
The latter inequality give us a motivation to consider $u^*_t$ as
the most likely path to absorbtion of the normed process $x^K_t$.

Note that calculations for $\P(\tau_0\le T)$ on logarithmic scale
requires a non-standard technique. The set $ \{\tau_0\le
T\}=\big\{(x^K_t)_{t\in[0,T]}\in\mathsf{D}\big\} $, where
$$
\mathsf{D}=\big\{u\in \mathbb{C}_{[0,T]}: u_0=1;
u_t=u_{\theta(u)\wedge t, \ \theta(u)=\inf\{t:u_t=0\}\le T}\big\}.
$$
$\mathsf{D}$ is closed in the uniform metric ($\varrho$) in the
space $\mathbb{C}_{[0,T]}$ of continuous functions on $[0,T]$.
Hence, the upper limit $
\varlimsup\limits_{K\to\infty}\frac{1}{K^{2(1-\gamma)}}\log\P\big(\tau_0\le
T\big) $ is done according to the LDP technique. However,
$\mathsf{D}$ has an empty interior. This fact prevents us to use
the LDP technique for the lower bound
$\varliminf\limits_{K\to\infty}\frac{1}{K^{2(1-\gamma)}}\log\P(\tau_0\le
T)$. Nevertheless, we obtain this lower  by using   an
  inclusion $ \{\tau_0\le T\} \supseteq
\Big\{\frac{1}{K^{1-\gamma}}M_T<-1 \Big\}, $ where $M_T$ is
defined in \eqref{Meq}. The probability of the  latter is  easily
computable, and gives a surprising result
$$
\lim_{K\to
\infty}\frac{1}{K^{2(1-\gamma)}}\log\Big\{\frac{1}{K^{1-\gamma}}M_T<-1
\Big\}=
\varlimsup\limits_{K\to\infty}\frac{1}{K^{2(1-\gamma)}}\log\P\big(\tau_0\le
T\big).
$$
Of course,
$$
\lim_{K\to
\infty}\frac{1}{K^{2(1-\gamma)}}\log\Big\{\frac{1}{K^{1-\gamma}}M_T<-1
\Big\}\le
\varliminf\limits_{K\to\infty}\frac{1}{K^{2(1-\gamma)}}\log\P(\tau_0\le
T)
$$
which together the above establish the desired limit. This
trick is of independent interest and might be useful for
establishing LDP in other problems.

\section{\bf Asymptotic of $\b{\mathsf{P}(\tau_0\le T)}$  as $\b{K\to\infty}$ on logarithmic scale}

The random process $M_t$ (see \eqref{Meq}) is a Gaussian
martingale with the variation process $\langle
M\rangle_t=\mathsf{E}M^2_t$:
\begin{equation}\label{eq:Ma}
\langle M\rangle_t=\int_0^t\sigma^2(1-\gamma)^2e^{-2(1-\gamma)\mu s}ds.
\end{equation}

\medskip
\noindent
\begin{theorem}\label{theo-2.1}
For any $T>0$,
$$
\lim\limits_{K\to\infty}
\frac{1}{K^{2(1-\gamma)}}\log\mathsf{P}(\tau_{0 }\le T)
=-\frac{1}{2\langle M\rangle_T}.
$$
\end{theorem}

\begin{proof}
To apply the It\^o formula   in a vicinity of $\tau_0$, let us
define a stopping time $ \tau_\varepsilon=\inf\{t\le
T:x^K_t=\varepsilon\}, \ \varepsilon>0. $ Now, by  It\^o's
formula, applied to $(x^K_t)^{1-\gamma}$,
$t\le\tau_\varepsilon\wedge T$, we find that
\begin{gather}
(x^K_t)^{1-\gamma}=1+\int_0^t(1-\gamma)\mu (x^K_s)^{1-\gamma}ds+
\int_0^t(1-\gamma)\sigma\frac{dB_s}{K^{1-\gamma}}
\nonumber\\
-\frac{1}{2}\int_0^t\gamma(1-\gamma)\frac{\sigma^2}{K^{2(1-\gamma)}}\frac{1}{(x^K_s)^{1-\gamma}}ds
\label{eq:2.4a}
\end{gather}
and, in turn,
\begin{gather*}
(x^K_{\tau_\varepsilon\wedge T})^{1-\gamma}_{\tau_\varepsilon\wedge T}e^{-(1-\gamma)\mu
(\tau_\varepsilon\wedge T)}
\\
 +\int_0^{\tau_\varepsilon\wedge T}
\frac{\sigma^2}{2K^{2(1-\gamma)}}\gamma(1-\gamma)\frac{e^{(1-\gamma)\mu s}}
{(x^K_s)^{1-\gamma}}ds=1+\frac{1}{K^{1-\gamma}}M_{\tau_\varepsilon\wedge T}.
\end{gather*}
In view of
$
\lim_{\varepsilon\to 0}M_{\tau_\varepsilon\wedge T}=M_{\tau_0\wedge T}
$
a.s. and the monotone convergence theorem
\begin{multline*}
\lim_{\varepsilon\to 0}\int_0^{\tau_\varepsilon\wedge T}
\frac{\sigma^2}{2K^{2(1-\gamma)}}\gamma(1-\gamma)\frac{e^{(1-\gamma)\mu s}}
{(x^K_s)^{1-\gamma}}ds
\\
=\int_{[0,\tau_0\wedge T)}
\frac{\sigma^2}{2K^{2(1-\gamma)}}\gamma(1-\gamma)\frac{e^{(1-\gamma)\mu s}}
{(x^K_s)^{1-\gamma}}ds, \ \text{a.s.},
\end{multline*}
in both sides of the above equality $\lim_{\varepsilon\to 0}$ is applicable, that is,
we have
\begin{gather}
0\le (x^K_{\tau_0\wedge T})^{1-\gamma}e^{-(1-\gamma)\mu
(\tau_0\wedge T)}
\nonumber\\
+\int_{[0,\tau_0\wedge T)}
\frac{\sigma^2}{2K^{2(1-\gamma)}}\gamma(1-\gamma)\frac{e^{(1-\gamma)\mu s}}
{(x^K)^{1-\gamma}_s}ds \ =1+\frac{1}{K^{1-\gamma}}M_{\tau_0\wedge T}.
\label{eq:2.1y}
\end{gather}
\eqref{eq:2.1y} implies
$
1+\frac{1}{K^{1-\gamma}}M_{\tau_0\wedge T}\ge 0.
$
If  $\omega\in \{\tau_0>T\}$, then
$1+\frac{1}{K^{1-\gamma}}M_T(\omega)\ge 0.$ In other words,
$\{\tau_0>T\}\subset \{1+\frac{1}{K^{1-\gamma}}M_T\ge 0\}$, and so
we obtain inclusion
\begin{equation}\label{eq:2.4b}
\{\tau_0\le T\} \supseteq \Big\{\frac{1}{K^{1-\gamma}}M_T+1< 0
\Big\}.
\end{equation}
It is well known that the families
$\big\{\frac{1}{K^{1-\gamma}}M_T\big\}_{K\to \infty}$ obeys LDP in the metric space $(\mathbb{R},\rho)$
($\rho$ is the Euclidian metric) with the rate speed $\frac{1}{K^{2(1-\gamma)}}$
and the rate function $ I(v)=\frac{v^2}{2\langle M\rangle_T}. $ In accordance with the large deviation
theory,
\begin{equation}\label{eq:lwb}
\varliminf_{K\to\infty}\frac{1}{K^{2(1-\gamma)}}\log\mathsf{P}\big(\tau\le T_0\big)\ge -\inf_{v:v+1\le 0}I(v)
=-\frac{1}{2\langle M\rangle_T}.
\end{equation}

A verification of the upper bound
\begin{equation}\label{eq:upb}
\varlimsup_{K\to\infty}\frac{1}{K^{2(1-\gamma)}}\log\mathsf{P}\big(\tau\le T_0\big)\le
-\frac{1}{2\langle M\rangle_T}
\end{equation}
is more involved. We select a set
$$
\mathsf{D}=\big\{u\in \mathbb{C}_{[0,T]}: u_0=1; u_t=u_{\theta(u)\wedge t, \
\theta(u)=\inf\{t:u_t=0\}\le T}\big\}
$$
which is closed in the uniform metric ($\varrho$) related to the
space  $\mathbb{C}_{[0,T]}$ of continuous functions on $[0,T]$.
Obviously, $ \{\tau_0\le
T\}\subseteq\big\{(x^K_t)_{t\in[0,T]}\in\mathsf{D}\big\} $, which
suggests to find
\begin{equation}\label{00ps}
\varlimsup_{K\to\infty}\frac{1}{K^{2(1-\gamma)}}\log\mathsf{P}
\Big((x^K_t)_{t\in[0,T]}\in\mathsf{D}\Big).
\end{equation}
The most convenient tool to this asymptotic analysis is LDP for family
$
\big\{\big(x^K_t\big)_{t\in[0,T]}\big\}_{K\to\infty}
$
having the speed rate
$
\frac{1}{K^{2(1-\gamma)}} (!)
$
and the rate function
\begin{equation*}
J_{T}(u)= \left\{
\begin{array}{lll}
\frac{1}{2\sigma^2}\displaystyle{\int_0^{\theta(u)\wedge T}}\Big(\frac{\dot{u}_t-\mu
u_t}{u^\gamma_t}\Big)^2dt, &\substack{ u_0=1\\du_t=\dot{u}_tdt, \\
\int\limits_{[0,\theta(u)\wedge T]}[\frac{\dot{u}_t-\mu u_t}{u^\gamma}]^2dt<\infty}
\\
\infty, & \text{otherwise}
\end{array}
\right.
\end{equation*}
(Theorem \ref{theo-A.1}). In accordance to the large deviation theory
\begin{equation*}
\varlimsup_{K\to\infty}\frac{1}{K^{2(1-\gamma)}}\log\P\Big((x^K_t)_{t\in[0,T]}\in\mathsf{D}\Big)
\le -\inf_{u\in\mathsf{D}}J_T(u),
\end{equation*}
so that, it remains to prove
$$
\inf_{u\in\mathsf{D}}J_T(u)=\frac{1}{2\langle M\rangle_T}.
$$
A minimization procedure of $J(u)$ in $u\in\mathsf{D}$ exclude from consideration functions
$u_t$ with $du_t\not\ll dt$ and
$
\int_0^{\theta(u)\wedge T}\big[\frac{\dot{u}_t-\mu u_t}{u^{\gamma}_t}\big]dt=\infty.
$
This minimization is realized with a help of a
specific deterministic control problem
with a control action $w_t$ and a controlled process
$u_t$, being the solution of differential equation
\begin{equation}\label{eq:nonl}
\dot{u}_t=\mu u_t+\sigma u^\gamma_tw_t, \ t\le \theta(u)\wedge T
\end{equation}
subject to the initial condition $u_0=1$. Obviously, the function $u_t$ belongs to $\mathsf{D}$.
The pair $(w^*_t,\theta(u^*))$, with $u^*(t)$ related to $w^*_t$, is said to be optimal if
$$
\int_0^{\theta(u^*)\wedge T}(w^*_t)^2dt\le \int_0^{\theta(u)\wedge
T} w^2_tdt
$$
for any pair $(w_t,u_t)$ with $ \int_0^{\theta(u)\wedge
T}w^2_tdt<\infty$. Technically, it is convenient to use the
following change of variables: $v_t=u^{1-\gamma}_t$ enables us to
reduce the problem to a linear differential equation
\begin{equation}\label{eq:lin}
\dot{v}_t=\mu(1-\gamma) v_t+\sigma(1-\gamma)w_t;
\end{equation}
instead of nonlinear one \eqref{eq:nonl}, subject to the initial condition $v_0=1$. We shall exploit
also the property
$
v_{\theta(u)}
  \begin{cases}
     =0, & u_{\theta(u)}=0 \\
    >0 , & u_{\theta(u)}>0.
  \end{cases}
$ The explicit  solution of equation \eqref{eq:lin}, under the
assumption $\theta(u)\le T$, implies:
$$
0=v_{\theta(u)}e^{-\mu(1-\gamma)t}=\Big[1+\sigma(1-\gamma)\int_{[0,\theta(u)\wedge
T)} e^{-\mu(1-\gamma) s}w_sds\Big]
$$
or, equivalently, the equality:
\begin{equation}\label{equ}
-\frac{1}{\sigma(1-\gamma)}=\int_{[0,
\theta(u)\wedge T)} e^{-\mu(1-\gamma)s}w_sds
\end{equation}
that, due to the Cauchy-Schwarz inequality, can be transformed into the inequality:
\begin{gather*}
\int_{[0, \theta(u)\wedge T)}
w^2_tdt\ge\frac{2\mu}{\sigma^2(1-\gamma)[1-e^{-2\mu(1-\gamma)\theta(u)}]}
\ge
\frac{2\mu}{\sigma^2(1-\gamma)[1-e^{-2\mu(1-\gamma) T}]}.
\end{gather*}
The choice of  $w^*_t$ is conditioned by two requirements:

1) \eqref{equ} remains valid for $w_t$ replaced by $w^*_t$

2)
$
\int_{[0, \theta(u^*)\wedge T)}(w^*_t)^2dt=
\frac{2\mu}{\sigma^2(1-\gamma)[1-e^{-2\mu(1-\gamma) T}]}.
$

\smallskip
\noindent
Both requirements are satisfied for
$
w^*_t=-\frac{1}{\sigma} \frac{2\mu}{1-e^{-2\mu(1-\gamma)T}}e^{-\mu(1-\gamma)t}.
$
Hence,
$
\frac{1}{2}\int_0^T(w^*_t)^2dt=\frac{1}{2\langle M\rangle_T}.
$
\end{proof}

\begin{remark}
$ u_t=-\frac{\langle M\rangle_t}{\langle M\rangle}_T. $
\end{remark}
\begin{remark}\label{rem-1} The fact that the   random variable $M_T$ is gaussian  with parameters
$(0,\langle M\rangle_t)$ and \eqref{eq:2.4b} yield for any $K>0$,
$ \P\big(\tau_0\le T\big)\ge \P\big(M_T\le-K^{1-\gamma}\big)  $.
\end{remark}

\section{\bf Most likely path to ruin of the normed process $\b{x^K_t}$}
Since $u^*_t\equiv (v^*_t)^{1/\gamma}$, where $v^*_t$ solves the
differential equation
$$
\dot{v}^*_t=\mu(1-\gamma) v^*_t+\sigma(1-\gamma)w^*_t
$$
with $v^*_0=1$, we find that
\begin{equation*}
u^*_t=
   e^{\mu t}\left[1-
   \frac{1-\frac{2\mu\langle M\rangle_t}{\sigma^2(1-\gamma)}}
   {1-\frac{2\mu\langle
   M\rangle_T}{\sigma^2(1-\gamma)}}\right]^{1/(1-\gamma)}\equiv
   e^{\mu t}\left[1-
   \frac{e^{-2(1-\gamma)\mu t}}
   {e^{-2(1-\gamma)\mu T}}\right]^{1/(1-\gamma)}.
\end{equation*}
On the other hand, in accordance with Theorem \ref{theo-A.1} for $u^*$, we have
\begin{equation*}
\lim_{\delta\to 0}\lim_{K\to\infty}\frac{1}{K^{2(1-\gamma)}}\log\P\bigg(\sup_{t\in[0,T]}|x^K_t-u^*_t|\le
\delta\bigg)
=-J_T(u^*).
\end{equation*}
At the same time for any $u\in\mathsf{D}$, Theorem \ref{theo-A.1} provides
\begin{multline*}
\lim_{\delta\to 0}\lim_{K\to\infty}\frac{1}{K^{2(1-\gamma)}}\log\P\bigg(\sup_{t\in[0,T]}|x^K_t-u_t|\le
\delta\bigg)
\\
=-J_T(u)\le-\inf\{u\in\mathsf{D}^\circ \le-J_T(u^*).
\end{multline*}

Consequently, the function $u^*_t$ can be considered as the most
likely path to ruin of the normed process $x^K_t$ on time interval
$[0,T]$.

\appendix
\section{\bf LDP for the family $\b{\{(x^K_t)_{t\in [0,T]}\}_{K\to\infty}}$}
\label{sec-5}

The family $\{(x^K_t)_{t\in [0,T]}\}_{K\to\infty}$ is in Freidlin-Wentzell's framework \cite{FW}.
In our setting, we take into account that the random process
$x^K_t$ is absorbed at the stopping time $\tau_0$, so that its paths belong to a subspace
$\mathbb{C}^{\rm abc}_{[0,T]}(\mathbb{R}_+)$ of
$\mathbb{C}_{[0,T]}(\mathbb{R}_+)$ the  space of continuous nonnegative
functions $u_t=u_{t\wedge\theta(u)}$, where
$\theta(u)=\inf\{t\le T:u_t=0\}$. The subspace
$\mathbb{C}^{\rm abc}_{[0,T]}(\mathbb{R}_+)$ is closed in the uniform metric $\varrho$ and,
consequently, it suffices to analyze the LDP in the metric space
$(\mathbb{C}^{\rm abc}_{[0,T]}(\mathbb{R}_+),\varrho)$. The use of
$(\mathbb{C}^{\rm abc}_{[0,T]}(\mathbb{R}_+),\varrho)$ instead of
$\mathbb{C}_{[0,T]}(\mathbb{R}_+)$ enables us to apply standard approach to LDP proof adding
a few simplest details only.

\begin{theorem}\label{theo-A.1}
The family $\big\{(x^K_t)_{t\ge 0}\}_{K\to \infty}$ obeys LDP in the metric space
$(\mathbb{C}^{\rm abc}_{[0,T]}(\mathbb{R}_+),\varrho)$ with the speed rate
$\frac{1}{K^{2(1-\gamma)}}$ and the rate function
\begin{equation*}
J_{T}(u)= \left\{
\begin{array}{lll}
\frac{1}{2\sigma^2}\displaystyle{\int_0^{\theta(u)\wedge T}}\Big(\frac{\dot{u}_t-\mu
u_t}{u^\gamma_t}\Big)^2dt, &\substack{ u_0=1\\du_t=\dot{u}_tdt, \\
\int\limits_{[0,\theta(u)\wedge T]}[\frac{\dot{u}_t-\mu u_t}{u^\gamma}]^2dt<\infty}
\\
\infty, & \text{otherwise}.
\end{array}
\right.
\end{equation*}
\end{theorem}

\begin{proof} The family $\big\{(x^K_t)_{t\ge 0}\}_{K\to \infty}$ is exponentially tight
(see, e.g., theorems 1.3 and 3.1, Liptser and Puhalskii, \cite{LP}), that is,
\begin{align}\label{b.1}
& \lim_{C\to\infty}\varlimsup_{K\to
0}\frac{1}{K^{2(1-\gamma)}}\log\mathsf{P} \Big(\sup_{t\in[0,T]} x^K_t\ge
C\Big)=-\infty,
\\
&\lim_{\Delta\to 0}\varlimsup_{K\to 0}\sup_{\vartheta\le
T}\frac{1}{K^{2(1-\gamma)}}\log\mathsf{P}\Big(\sup_{t\in[0,\triangle]}
|x^K_{\vartheta+t}-x^K_\vartheta|\ge
\eta\Big) =-\infty,  \label{b.2}
\end{align}
where $\eta$ is arbitrary number and $\vartheta$ is stopping time
relative to a corresponding filtration. In \eqref{b.1}, without
loss generality  $x^K_t$ might be replaced by $(x^K_t)^{1-\gamma}$
what makes possible, in accordance with \eqref{eq:2.4a}, to use
the inequality $ (x^K_t)^{1-\gamma}\le e^{(1-\gamma)\mu t}\big[1+
\int_0^te^{(-1-\gamma)\mu
s}(1-\gamma)\sigma\frac{dB_s}{K^{1-\gamma}}\big], $ making the
proof transparent. Due to \eqref{b.1}, the condition from
\eqref{b.2} can be replaced by an easy provable condition (here $
\frak{A}_C=\big\{\sup_{t\le T}x^\varepsilon_t\le C\big\}$):
$$
\lim_{\triangle\to 0}\varlimsup_{K\to 0}\sup_{\vartheta\le
T}\frac{1}{K^{2(1-\gamma)}}\log\mathsf{P}\Big(\sup_{[0,\triangle]}|x^K_{\vartheta+t}-x^K_\vartheta|\ge
\eta, \frak{A}_C\Big) =-\infty, \ \forall \ C>0.
$$
For $\theta(u)>T$, the proof of local LDP:
\begin{equation*}
\lim_{\delta\to
0}\lim_{K\to\infty}\frac{1}{K^{2(1-\gamma)}}\log\mathsf{P}\Big(\sup_{t\in[0,T]}
\big|x^K_t-u_t\big|\le\delta\Big)=-J_T(u)
\end{equation*}
does not different from standard one and is omitted.
The case of $J_T(u)=-\infty$, including $u_0\ne 1$, $du_t\not\ll dt$,
is analyzed in a standard way and is omitted too.

\medskip
The analysis of
``$u_0=1$, $du_t=\dot{u}_tdt$,
$
\int_0^{\theta(u)\wedge T}\big(\frac{\dot{u}_s-\mu
u_s}{u^\gamma_s}\big)^2ds<\infty,$ $\theta(u)\le T$'' is based on the following result.
\begin{proposition}\label{pro-2}{\rm [Dupuis, Ellis \cite{DupEll}, Теорема A.6.3.]}
For any absolutely continuous function
$u=(u_t)_{t\in[0,T]}$, mapping  $[0,T]$ into $\mathbb{R}$, and any $a\in \mathbb{R}$
$$
\int_0^TI_{\big\{\substack{u_t=a\\ \dot{u}_t\neq 0}\big\}}dt=0.
$$
\end{proposition}

{\bf Local LDP upper bound.}
Set $u^n_t=\frac{1}{n}\vee u_t$ and notice that $\theta(u^n)>T$.
Moreover, $u^n_0=1$, $du^n_t=\dot{u}^n_tdt$ and, due to Proposition
\ref{pro-2},
$
\dot{u}^n_t=\dot{u}_tI_{\{u_t>\frac{1}{n}\}}ds
$
and also
$
\int_0^{\theta(u^n)\wedge T}\big(\frac{\dot{u}^n_s-\mu
u^n_s}{(u^n)^\gamma_s}\big)^2ds<\infty.
$

Since $\tau^n=\inf\{t:u_t\le \frac{1}{n}\}\to \theta(u)$, $n\to\infty$, we find that
\begin{align*}
&\varlimsup_{\delta\to
0}\varlimsup_{K\to\infty}\frac{1}{K^{2(1-\gamma)}}\log\mathsf{P}\Big(\sup_{t\in[0,T]}
\big|x^K_t-u_t\big|\le\delta\Big)
\\
&\le \varlimsup_{\delta\to
0}\varlimsup_{K\to\infty}\frac{1}{K^{2(1-\gamma)}}\log\mathsf{P}\Big(\sup_{t\in[0,\tau^n\wedge T]}
\big|x^K_t-u_t\big|\le\delta\Big)
\\
&\le-\frac{1}{2\sigma^2}\int_0^{\tau^n\wedge T}\Big(\frac{\dot{u}^n_t-\mu u^n_t}
{(u^n_t)^\gamma}\Big)^2dt=-\frac{1}{2\sigma^2}\int_0^{\tau^n\wedge T}
\Big(\frac{\dot{u}_t-\mu u_t}{(u_t)^\gamma}
\Big)^2dt
\\
&\xrightarrow[n\to\infty]{}-\frac{1}{2}\int_0^{\theta(u)\wedge T}\Big(\frac{\dot{u}_t-\mu u_t}{(u_t)^\gamma}
\Big)^2dt.
\end{align*}

{\bf Local LDP lower bound.}
With $\phi>\delta>0$, write
\begin{align*}
&\Big\{\sup_{t\in[0,T]}|x^K_t-u^n_t|\le \delta\Big\}
\\
&= \Big\{\sup_{t\in[0,T]}|x^K_t-u^n_t|\le
\delta\Big\}\cap\Big\{\sup_{t\in[0,T]}|u^n_t-u_t|\le \phi\Big\}
\\
&\quad \bigcup\Big\{\sup_{t\in[0,T]}|x^K_t-u^n_t|\le
\delta\Big\}\cap\Big\{\sup_{t\in[0,T]}|u^n_t-u_t|>\phi\Big\}
\\
&\subseteq
\Big\{\sup_{t\in[0,T]}|x^K_t-u_t|\le\phi+\delta\Big\}
\cap\Big\{\sup_{t\in[0,T]}|u^n_t-u_t|>\phi\Big\}
\\
&\subseteq \Big\{\sup_{t\in[0,T]}|x^K_t-u_t|\le
2\phi\Big\} \cup\Big\{\sup_{t\in[0,T]}|u^n_t-u_t|> \phi\Big\}.
\end{align*}
For fixed $\phi$, there exits a number $n_\phi>\frac{1}{\phi}$ such that for any $n\ge n_\phi$
the set
$
\big\{\sup_{t\in[0,T]}|u^n_t-u_t|> \phi\big\}=\varnothing.
$
Therefore, for sufficiently large numbers $n$,
\begin{equation*}
\mathsf{P}\Big(\sup_{t\in[0,T]}|x^K_t-u_t|\le
2\phi\Big)
\\
\ge
\mathsf{P}\Big(\sup_{t\in[0,T]}|x^K_t-u^n_t|\le
\delta\Big).
\end{equation*}
Hence and by Proposition \ref{pro-2}, a chain of lower bounds holds,
 \begin{align*}
&\varliminf_{K\to\infty}\frac{1}{K^{2(1-\gamma)}}\log\mathsf{P}\Big(\sup_{t\in[0,T]}|x^K_t-u_t|\le
2\phi\Big)
\\
&\ge\varliminf_{\delta\to 0}\varliminf_{K\to\infty}\frac{1}{K^{2(1-\gamma}}\log\mathsf{P}
\Big(\sup_{t\in[0,T]}
|x^K_t-u^n_t|\le \delta\Big)
\\
&\ge-\frac{1}{2\sigma^2}\int_0^{\tau^n\wedge T} \Big(\frac{\dot{u}_s-\mu
u_s}{u^\gamma_s}\Big)^2ds-\frac{1}{2\sigma^2}\int_{\tau^n\wedge T}^T\frac{\mu^2}{n^{2(1-\gamma)}}
\\
&\xrightarrow[n\to\infty]{}-\frac{1}{2\sigma^2}\int_0^{\theta(u)\wedge T} \Big(\frac{\dot{u}_s-\mu
u_s}{u^\gamma_s}\Big)^2ds,
\end{align*}
providing
\begin{multline*}
\varlimsup_{\phi\to 0}\varlimsup_{K\to\infty}
\frac{1}{K^{2(1-\gamma)}}\log \mathsf{P}\Big(\sup_{t\in[0,T]}|x^K_t-u_t|\le 2\phi\Big)
\\
 \le -\frac{1}{2\sigma^2}
\int_0^{\theta(u)\wedge T}\Big(\frac{\dot{u}_s-\mu
u_s}{u^\gamma_s}\Big)^2ds.
\end{multline*}
\end{proof}

\end{document}